\documentclass[twoside]{article}

\usepackage[accepted]{aistats2019}
\usepackage{tikz}
\usetikzlibrary {shapes}
\usetikzlibrary {arrows}
\usetikzlibrary {positioning}
\usetikzlibrary {calc}
\usetikzlibrary{fit}

\usepackage{graphicx}

\usepackage{here}

\usepackage{latexsym}

\usepackage{amsbsy}

\usepackage{amsmath,amsthm,amssymb, amsfonts}

\usepackage{dsfont}

\newcommand{\dse}{\,\mbox{$\perp$}\,}

\newcommand{\cip}{\mbox{\,$\perp\!\!\!\perp$\,}}
\newcommand{\sk}{\mathrm{sk}}

\newcommand{\cd}{\,|\,}

\newcommand{\mtp}{{\rm MTP}_{2}}
\newcommand{\ci}{\mbox{\protect $\: \perp \hspace{-2.3ex}
\perp$ }}

\newcommand{\notci}{\nolinebreak{\not\hspace{-1.5mm}\ci}}

\newcommand{\nn}[0]{\hspace*{.7em}}

\newcommand{\node}{\mbox {\LARGE
{$\mbox{$\circ$}$}}}

\newcommand{\ful}{\mbox{$\, \frac{ \nn \nn \;}{ \nn \nn
}$}}

\newcommand{\fra}{\mbox{$\hspace{.05em} \frac{\nn
\nn}{\nn
}\!\!\!\!\! \succ \! \hspace{.25ex}$}}

\newcommand{\arc}{\mbox{$\hspace{.06em} \prec
\!\!\!\!\!\frac{\nn \nn}{\nn}
\!\!\!\!\!
\succ\! \hspace{.25ex}$}}

\newtheorem{prop}{Proposition}
\newtheorem{coro}{Corollary}

\newtheorem{lemma}{Lemma}

\newtheorem{theorem}{Theorem}
\newtheorem{example}{Example}

\usepackage[round]{natbib}

\begin{document}

%

%

\twocolumn[

\aistatstitle{Markov Properties of Discrete Determinantal Point Processes}

\aistatsauthor{Kayvan Sadeghi  \And Alessandro Rinaldo }

\aistatsaddress{ University College London \And  Carnegie Mellon University } ]

\begin{abstract}
Determinantal point processes (DPPs) are probabilistic models for repulsion. When used to represent the occurrence of random subsets of a finite base set, DPPs allow to model global negative associations in a mathematically elegant and direct way. Discrete DPPs have become popular and computationally tractable models for solving several machine learning tasks that require the selection of diverse objects, and have been successfully applied in numerous real-life problems. Despite their popularity, the statistical properties of such models have not been adequately explored. In this note, we derive the Markov properties of discrete DPPs and show how they can be expressed using graphical models.
\end{abstract}

\section{Introduction}

Determinantal point processes (DPPs) are stochastic processes for repulsion that model in a mathematically elegant and general way negative association. Originally developed in quantum physics, DPPs arise naturally in other areas, such as combinatorics, random matrix theory, probability and algebra; see, e.g.,  \cite{lyo03}, \cite{bor09}, \cite{borodin}, \cite{hough2006} \cite{Johanssonrandommatrices} and references therein.
In particular, DPPs on finite base sets  model the emergence of random subsets comprised by objects with a higher degree of mutual diversity that would  have been  otherwise observed through independent selections.  In detail, for a given finite set $V$ and a positive definite matrix $K$ with eigenvalues bounded by $1$, a discrete DPP $Y$ is a random subset of $V$ such that, for any $A \subset V$,
\[
\mathbb{P}\left( A \subseteq  Y \right) = \det(K_A),
\]
where $K_A$ is the principal minor of $K$ corresponding to the coordinates in $A$ and $\det(K_{\emptyset}) = 1$.

In their seminal work, \cite{kul12} have demonstrated the properties, expressive power, and computational ease of discrete DPPs for machine learning tasks aimed at producing diverse subset selections, and have exemplified their use in several real-life problems. Discrete DPPs have now been successfully employed in numerous areas, ranging from document summarization, signal processing, neuroscience, image segmentation, spatial statistics, survey sampling, and wireless network; see, e.g., \cite{bru17} and references therein for a list of applications of discrete DPPs. Despite their popularity, the formal statistical properties of such models and their statistical inference remain largely unexplored; see \cite{bru17} and \cite{urschel} for recent works on maximum likelihood estimation, sample complexity guarantees and efficient algorithms in these models.

The main goal of this paper is to characterize the independence properties of discrete DPPs, and to express them trough the formalism of graphical modeling. By representing DPPs over finite sets as random binary vectors, we are able to easily derive the Markov properties of the corresponding distributions, using both bidirected and undirected graphical models. In particular, we make the following contributions:
\begin{enumerate}
\item We show that the distributions of discrete DPPs are Markov with respect to bidirected graphical models with edges corresponding to the non-zero entries of $K$. Such models are curved exponential families \citep[see, e.g.]{drt08} and their Markov properties as well as other statistical properties are well-studied; see \cite{kau96,rov13}.
\item For a {\it conditional} discrete DPPs $Y$, which is a DPP over a fixed non-empty subset $A$ of $V$ conditionally  on $A^c \subset Y$, we further show that the corresponding conditional distribution is Markov with respect to an undirected graph whose edges instead correspond to the non-zero entries of $K^{-1}$. In particular, such conditional models are linear exponential families. This particular type of independence properties are known as \emph{context-specific} independence models, since they only apply to a distinguished subset of conditional distributions.

\item We further explore conditions that guarantee faithfulness in either one of the above cases.
\end{enumerate}
Interestingly, we have found that DPP graphical models act similarly to the Gaussian graphical models, whereby the kernel matrix $K$ plays a role analogous to that of the covariance matrix. Indeed, both $K$ and $K^{-1}$ express Markov properties that can be easily elicited from the corresponding bidirected and undirected graphs.
To the best of our knowledge, results of this form are novel. Previous work on this topic include the results  by \cite{tad14}, which, among the other things, relate the conditional distribution of  a discrete DPPs $Y$  given $C \subset Y $ to the Schur complement of $\bar{K}$ in $K$.

\section{Preliminaries}
We first provide some background on determinantal point process on finite sets and review basic concepts and some recent results from the theory of probabilistic independence models and graphical models. 

\subsection{Determinantal point processes}
Let $V$ be a finite set of cardinality $n$, which, without loss of generality, may be assumed to coincide with $\{1,\ldots,n\}$. A \emph{(discrete) determinantal point process (DPP)} $Y$ on $V$ is random variables taking values in  $2^{V}$, the set of all subsets of $V$, such that
\begin{equation}\label{eq:first}
\mathbb{P}\left( A \subseteq  Y \right) = \det(K_A), \quad A \subseteq V, A \neq \emptyset,
\end{equation}

where $K \in \mathbb{R}^{V \times V}$ is a positive definite  matrix called the \emph{marginal kernel} or the simply \emph{kernel} of the process. Here, for any  non empty subset $A$ of $V$, $K_A$ denotes the restriction of $K$ to the entries $A \times A$, whereby we set $\det(K_{\emptyset}) = 1$. Similarly, for any non-empty subsets $A$ and $B$ of $V$, $K_{A,B}$ is  the restriction of $K$ to the row entries indexed by $A$ and column entries indexed by $B$.
By \eqref{eq:first}, all principal minors of $K$ must be nonnegative and no larger than $1$.  It is also easy to show  that the eigenvalues of $K$ are bounded above by $1$.  In fact, these requirements are also sufficient:  any matrix $K \in \mathbb{R}^{V \times V}$  such that  $0 \preceq K\preceq I$, where $I$ is the identity matrix of the appropriate size, defines a DPP  \citep[see, e.g., Theorem 2.3 in][]{kul12}.

For any pair $i$ and $j$ in $A$, (\ref{eq:first}) implies that
\[
\mathbb{P}( \{i,j\} \in Y) = \det(K_{\{i,j\}}) = K_{\{i\}} K_{\{j\}} - K^2_{\{i,j\}},
\]
so that $\mathbb{P}( \{i,j\} \in Y) \leq \mathbb{P}( i \in Y) \mathbb{P}( j \in Y)$.
\subsubsection{Marginalization, Conditioning and Complements} The class of probability distributions of DPPs over $V$ is closed under marginalization, conditioning and complements. The first property is immediate: if $Y$ is a DPP over $V$ with kernel $K$ and $A \subset V$ is non-empty, then $Y_A := Y \cap A$ is a DPP over $A$ with kernel $K_A$. As for the conditioning, \cite{tad14} has showed that for any non-empty $A \subset V$, the conditional distribution of $Y_{\bar{A}}$ given that $Y_A = 1_A$, where $\bar{A} = V \setminus A$, is that of a DPP over $\bar{A}$ with kernel given by $K^A:=K_{\bar{A}}-K_{\bar{A},A}K^{-1}_AK_{A,\bar{A}}$, the \emph{Schur complement} of $K_A$ in $K$. In addition, the \emph{complement} $\bar{Y}$ of $Y$ given by $\bar{Y} = V \setminus Y$ is a DPP over $V$ with kernel $I - K$ \citep[see, e.g.][]{kul12}. Thats is,
\begin{equation}\label{eq:complement}
\mathbb{P}(A \subseteq \bar{Y}) = \det(I - K_A), \quad A \subseteq V.
\end{equation}

\subsubsection{DPPs as Random Binary Vectors}
We will equivalently express the DPP $Y$ on $V$ as the $n$-dimensional binary random random vector $X_V = X = (X_1,\ldots,X_n)$ taking values on $\{ 0,1 \}^{V}$ such that  $X_i = 1$ if and only if $i \in Y$.
Letting, for any non-empty $ A \subset V$, $X_A$ be the restriction of $X$ to the coordinates in $A$, the condition \eqref{eq:first} can be re-written as
\begin{equation}\label{eq:1}
q_A:=\mathbb{P}(X_A=1_A) = \det(K_A), \quad A \subseteq V, A \neq \emptyset,
\end{equation}
where $1_A$ is a vector with all entries equal to $1$ and indexed by $A$ and we set $q_{\emptyset} = 1$. Similarly, \eqref{eq:complement} becomes
\[
\mathbb{P}( X_A = 0_A) = \det(I - K_A).
\]
The parameters  $\{ q_A, A \subseteq V\}$ are easily seen to be marginal probabilities and, borrowing the terminology from \cite{drt08}, we will refer to them as the \emph{M\"{o}bius parameter}.
The joint probabilities of a finite DPP can be expressed in closed form using the so-called \emph{L-ensembles} (see \citealt{kul12}), assuming that the spectrum of $K$ does not contain $1$. In fact, they can also be obtained from the M\"{o}bius parameters using the inclusion-exclusion principle as follows:
$$P(X_A=1,X_{V\setminus A}=0)=\sum_{B \subseteq V:A\subseteq B}(-1)^{|B\setminus A|}q_B,$$
for all $A \subset V$. Notice that the above expression holds true even of the largest eigenvalue of $K$ is $1$.

As remarked above, the set of distributions over $\{0,1\}^V$ that correspond to DPPs $X_V$ over $V$ is closed under marginalization, conditioning  on  events of the form $X_A = 1_A$, for any non-empty $A \subset V$,  and complements, in the sense that $1_V - X_V$ also represents a DPP over $V$. The above operations can be combined to generated other DPPs.  In particular, this shows that DPPs are closed under conditioning on $X_A = 0_A$, the $|A|$-dimensional vector of all $0$'s, and consequently, since conditioning is commutative and transitive, DPPs are closed under conditioning on any event of the form $X_A = x_A$, for any non-empty $A \subset V$ and arbitrary $x_A \in \{0,1\}^A$.



\subsection{Independence models}\label{sec:propind}
An \emph{independence model} $\mathcal{J}$ over a finite set $V$ is a set of triples $\langle X,Y\cd Z\rangle$ (called \emph{independence statements}), where $X$, $Y$, and $Z$ are disjoint subsets of $V$; $Z$
may be empty, and $\langle \varnothing,Y\cd Z\rangle$ and $\langle X,\varnothing\cd Z\rangle$ are always included in $\mathcal{J}$. The independence statement $\langle X,Y\cd Z\rangle$ is read as ``$X$ is independent of $Y$ given $Z$''. Independence models may in  general have a  probabilistic interpretation, but not necessarily. Similarly, not all independence models can be easily represented by graphs. For further discussion on general independence models, see \citet{stu05}.

\subsubsection{Probabilistic independence models}
In order to define probabilistic independence models, consider a set $V$ and a collection of random variables
$\{X_v \}_{v \in V}$ with state spaces $\{ \mathcal{X}_v \}_{v \in V}$ and joint distribution $P$. We let $X_A=\{X_v\}_{v\in A}$ for each non-empty subset $A$ of $V$. For disjoint subsets $A$, $B$, and $C$ of $V$
we use the short notation $A\cip B\cd C$ to denote that $X_A$ is \emph{conditionally independent of $X_B$ given $X_C$} \citep{daw79,lau96}, i.e.\ that for any measurable $\Omega\subseteq \mathcal{X}_A$ and $P$-almost all $x_B$ and $x_C$,
$$P(X_A \in \Omega\cd X_B=x_B, X_C=x_C)=P(X_A \in \Omega\cd X_C=x_C).$$
We can now induce an independence model $\mathcal{J}(P)$ by letting
\begin{displaymath}
\langle A,B\cd C\rangle\in \mathcal{J}(P) \text{ if and only if } A\ci B\cd C \text{ w.r.t.\ $P$}.
\end{displaymath}

If $A$, $B$, or $C$ has only one member $\{i\}$, $\{j\}$, or $\{k\}$, for better readability, we write $i\ci j\cd k$. We also write $A\ci B$ when $C=\varnothing$, which denotes the \emph{marginal independence} of $A$ and $B$.

A probabilistic independence model $\mathcal{J}(P)$ over a set $V$ is always a \emph{semi-graphoid} \citep{pea88}, i.e., it satisfies the four following properties for disjoint subsets $A$, $B$, $C$, and $D$ of $V$, which, as we shall see, makes them useful for graphical modeling:
 \begin{enumerate}
    \item $A\ci B\cd C$ if and only if $B\ci A\cd C$ (\emph{symmetry});
    \item if $A\ci B\cup D\cd C$ then $A\ci B\cd C$ and $A\ci D\cd C$ (\emph{decomposition});
    \item if $A\ci B\cup D\cd C$ then $A\ci B\cd C\cup D$ and $A\ci D\cd C\cup B$ (\emph{weak union});
    \item if $A\ci B\cd C\cup D$ and $A\ci D\cd C$ then $A\ci B\cup D\cd C$ (\emph{contraction}).
 \end{enumerate}
Notice that the reverse implication of contraction clearly holds by decomposition and weak union. A semi-graphoid for which the reverse implication of the weak union property holds is said to be a \emph{graphoid}; that is, it also satisfies
\begin{itemize}
	\item[5.] if $A\ci B\cd C\cup D$ and $A\ci D\cd C\cup B$ then $A\ci B\cup D\cd C$ (\emph{intersection}).
\end{itemize}
Furthermore, a graphoid or semi-graphoid for which the reverse implication of the decomposition property holds is said to be \emph{compositional}, that is, it also satisfies
\begin{itemize}
	\item[6.] if $A\ci B\cd C$ and $A\ci D\cd C$ then $A\ci B\cup D\cd C$ (\emph{composition}).
\end{itemize}
If, for example, $P$ has strictly positive density, the induced probabilistic independence model is always a graphoid; see e.g.\ Proposition 3.1 in \citet{lau96}. See also \citet{pet15} for a necessary and sufficient condition for $P$ in order for the intersection property to hold. If the distribution $P$ is a regular multivariate Gaussian distribution, $\mathcal{J}(P)$ is a compositional graphoid; e.g.\ see \citet{stu05}.
Probabilistic independence models with positive densities are not in general compositional; this only holds for special types of multivariate distributions such as, for example,  Gaussian distributions and the symmetric binary distributions used in \citet{wer09}.

Another important property that is not necessarily satisfied by probabilistic independence models, is \emph{singleton-transitivity} (also called \emph{weak transitivity} in \cite{pea88}, where it is shown that for Gaussian and binary distributions $P$, $\mathcal{J}(P)$ always satisfies it). For $i$, $j$, and $k$, single elements in $V$, and $C \subset V$,
 \begin{itemize}
	\item[7.] if $i\ci j\cd C$ and $i\ci j\cd C\cup\{k\}$ then $i\ci k\cd C$ or $j\ci k\cd C$ (singleton-transitivity).
\end{itemize}
In addition, we have the two following properties:
 \begin{itemize}
\item[8.] if $i\ci j\cd C$ then $i\ci j\cd C\cup\{k\}$ for every $k\in V\setminus \{i,j\}$ (\emph{upward-stability});
\item[9.] if $i\ci j\cd C$ then $i\ci j\cd C\setminus\{k\}$ for every $k\in V\setminus \{i,j\}$ (\emph{downward-stability}).
\end{itemize}

\subsubsection{Context-specific independence models}
In some situations, conditional independence may only hold when some variables are held fixed, an instance known as \emph{context-specific independence} \citep{bou96}. More precisely, let $A,B,D$, and $C$ be four disjoint sets of variables. $A$ and $B$ are conditionally independent given $D$ in context $X_C=x_c$ if, for any measurable $\Omega \subseteq \mathcal{X}_A$,
\[
P(X_A \in \Omega\cd X_B=x_B, X_C=x_C, X_D = x_D)\]
is equal to
\[P(X_A \in \Omega\cd X_C=x_C, X_D = x_D),
\]
 for $P$-almost all $x_B$ and $x_D$ such that $P( X_B=x_B, X_C=x_C, X_D = x_D) > 0$. Below we will express this property by writing $A \ci B \cd D \cup (C = x_C)$.
 When $D$ is empty, one simply says that $A$ and $B$ are independent in context $C=x_C$.

Now denote by $\mathcal{J}_{\cdot\cd X_C=x_C}(P)$ the \emph{context-specific independence model} induced by $P$, i.e., the set of all context-specific independence statements above. The following result on proving that certain axioms are satisfied by context-specific independence models is a key property of context-specific independence models. See also \citet{cor16}, for a more general discussion on the topic (although the following results are not explicitly proven in it).

Notice that the conditioning sets $C$ in the conditional independence statements appearing in the axioms for context-specific independence models $\mathcal{J}_{\cdot\cd X_{C'}=x_{C'}}$ could be of form $C'\cup C''$, where we condition on the context-specific $X_{C'}=x_{c'}$ and on the general $C''$. For example, the weak union property is defined as follows (other axioms are translated from the original definitions in a similar fashion):
 \begin{itemize}
    \item if $A\ci B\cup D\cd (C'=x_{C'})\cup C''$ then $A\ci B\cd (C'=x_{C'})\cup C''\cup D$ and $A\ci D\cd (C'=x_{C'})\cup C''\cup B$ (\emph{weak union}).
 \end{itemize}

\begin{lemma}\label{lem:2}
 A context-specific probabilistic independence model $\mathcal{J}_{\cdot\cd X_C=x_c}(P)$ satisfies the symmetry, decomposition, and weak union properties.
\end{lemma}
\begin{proof}
The proof is similar to that in the general case (see, e.g., \citealt{daw79}), where the value of the variables in conditioning set $C$ are set to a fixed $x_C$ (and in the weak-union case the values of the variables in $D$ are also set to a fix value $x_D$).
\end{proof}

However, the contraction property is not necessarily satisfied for context-specific independence models even in the case of binary random variables:
\begin{example}
Suppose that $A,B,C,D$ are four binary random variables. Let $P(C=0)=0$ (i.e., $C$ is deterministic and equal to $1$ almost surely). For any $x \in \{0,1\}^3$, set $P(x) = P(A=x_1,B=x_2,D=x_3)$. Next, define the probabilities
\begin{equation*}
\begin{split}
p_1=P(0,0,0), p_2=P(0,0,1), \\p_3=P(0,1,0), p_4=P(1,0,0),\\
p_5=P(0,1,1), p_6=P(1,0,1),\\ p_7=P(1,1,0), p_8=P(1,1,1).
\end{split}
 \end{equation*}
Now let $p_2=p_5=p_6=p_8=1/8$ and $p_1+p_3=p_4+p_7=1/4$, but let $p_1,p_3,p_4,p_7$ be all distinct. It is now easy to check that $A\ci B\cd D\cup C=(1,1)$ and $A\ci D\cd C=1$, but $A\notci B\cd C=1$.
\end{example}
\subsection{Graphical models}
Graphical models (see, e.g., \citealt{lau96}) are statistical models expressing
conditional independence clauses among a collection of random variables $X_V = (X_v, v
\in V)$ indexed by a finite set $V$. A graphical model is determined by a graph
$G=(V,E)$ over the
indexing set $V$, and the edge set $E$ (which may include edges of undirected, directed
or bidirected type) encodes conditional independence
relations among the variables, or {\it Markov properties}.

Using the formulation in (\ref{eq:1}), it is natural to define graphical models with node sets $X_1,\dots,X_n$ for DPPs.

\subsubsection{Undirected and bidirected graphical models}\label{sec:23}

For non-empty subsets $A$ and $B$ of $V$ and a (possibly empty) subset $S$ of
$V$, with $A$, $B$ and $S$ disjoint, we use the notation $A\ci B\cd S$ as shorthand for the conditional
independence relation $X_A\ci X_B\cd X_S$, thus identifying random variables
with their labels. When $S = \emptyset$ the independence relation is
intended as marginal independence.

For an undirected graph $G$, where all edges are undirected $\ful$, the  \emph{(global) Markov property} expresses that $A\ci B\cd S$ when every path between $A$ and $B$ has a vertex in $S$ or, in other words, $S$ \emph{separates} $A$ from $B$ in $G$. The \emph{pairwise Markov property} expresses that if two nodes $i$ and $j$ are non-adjacent, i.e.\ there is no edge between them, then $i\ci j\cd V\setminus \{i,j\}$. It is known that these two conditions are equivalent when the intersection property is satisfied; see \citet{pea88,lau96}.

For a \emph{directed acyclic graph} $G$, where all edges are directed $\fra$ and
there are no directed cycles, the  (global) Markov property expresses that
$A\ci B\cd S$ when  $A$ and $B$ are \emph{$d$-separated} \citep{pea88} by $S$; we
refrain from discussing the notion of $d$-separation here and refer to
\citet{lau96} for details.

%

We shall be interested in graphical models given by a \emph{bidirected graph} $G$ where all edges are  bidirected $\arc$. For such graphs the (global) Markov property \citep{cox93,kau96}  expresses that
$A\ci B\cd S$ when every path between $A$ and $B$ has a vertex outside $S\cup A\cup B$, i.e.\ $V\setminus(A\cup B\cup S)$ separates $A$ from $B$. Note the obvious duality between this and the Markov property for undirected graphs. The \emph{pairwise Markov property} expresses that if $i$ and $j$ are non-adjacent then $i\ci j$. It is known that these two conditions are equivalent for distributions that satisfy the composition property; see \citet{sadl14}.

\begin{figure}[htb]
\centering
\begin{tikzpicture}[node distance = 8mm and 8mm, minimum width = 6mm]
    \begin{scope}
      \tikzstyle{every node} = [shape = circle,
      font = \scriptsize,
      minimum height = 6mm,
      inner sep = 2pt,
      draw = black,
      fill = white,
      anchor = center],
      text centered]
      \node(i) at (0,0) {$i$};
      \node(j) [right = of i] {$j$};
      \node(l) [above = of j] {$l$};
			\node(k) [right = of j] {$k$};
      \node(i1) [right = 12mm of k] {$i$};
      \node(j1) [right = of i1] {$j$};
      \node(l1) [above = of j1] {$l$};
			\node(k1) [right = of j1] {$k$};
    \end{scope}
		
    \begin{scope}
    \tikzstyle{every node} = [node distance = 6mm and 6mm, minimum width = 6mm,
    font= \scriptsize,
      anchor = center,
      text centered]
\node(a) [below = 2mm  of j]{(a)};
\node(b) [below = 2mm  of j1]{(b)};

\end{scope}
    \begin{scope}
      \draw (i) -- (j);
      \draw (j) -- (l);
			\draw (j) -- (k);
    \end{scope}
    \begin{scope}[<->, > = latex']
    \draw (i1) -- (j1);
      \draw (j1) -- (l1);
			\draw (j1) -- (k1);
    \end{scope}

    \end{tikzpicture}
		\caption{{\footnotesize (a) An undirected graph. (b) a bidirected graph.}}\label{fig:001}
		\end{figure}
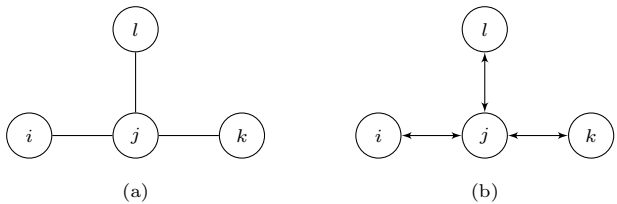

For example, in the undirected graph of Figure \ref{fig:001}(a),  the pairwise Markov property implies that $i\ci k\cd \{j,l\}$ and
the global Markov property implies that $\{i,l\}\ci k \cd j$ , whereas in the bidirected graph of Figure \ref{fig:001}(b), the pairwise Markov property implies that $i\ci k $ and the global Markov property implies that $\{i,l\}\ci k$.

Another Markov property for bidirected graphs is the so-called \emph{connected set Markov property} of \citet{drt08}, which is equivalent to the global Markov property. Lemma 1 in \citet{drt08} states that for discrete random variables taking values in the set $\mathcal{I}$, the connected set Markov property is satisfied if and only if for every disconnected set $D\subseteq V$ (i.e.\ for $D$ inducing a disconnected subgraph) such that $D = C_1\cup C_2\cup\cdots\cup C_r$, where $C_m$, $1\leq m\leq r$, are inclusion maximal connected sets, it holds that
\begin{equation}\label{eq:2}
\begin{split}
P(X_D = i_D) =  & P(X_{C_1} = i_{C_1}) \\
\times & P(X_{C_2} = i_{C_2})\cdots P(X_{C_r} = i_{C_r}),
\end{split}
\end{equation}
for every $i\in\mathcal{I}$.

If, for $P$ and undirected $G$, $A\dse_u B\cd C\iff A\ci B\cd C$ then we say that $P$ and $G$ are \emph{faithful}; Similarly, if, for $P$ and bidirected $G$, $A\dse_b B\cd C\iff A\ci B\cd C$ then $P$ and $G$ are \emph{faithful}. Hence, faithfulness implies Markovness, but not the other way around.

For a given probability distribution $P$, we define the \emph{skeleton} of
$P$, denoted by $\sk(P)$, to be the UG with node set $V$ such
that nodes $i$ and $j$ are not adjacent if and only if there is \emph{some} subset $C$
of $V$ so that $i\cip j\cd C$.  Thus, if $P$ is Markov with respect to an
undirected graph $G$ then $\sk(P)$ would be a subgraph of $G$ (by considering all remaining nodes to be in $C$); whereas if $P$ is
Markov with respect to a bidirected graph $G$ then $\sk(P)$ is a subgraph of $\sk(G)$ (by considering $C=\varnothing$).


In general a graph $G(P)$ is induced by $P$ with skeleton $\sk(P)$. For UGs, let $G_u(P)=\sk(P)$, whereas for BGs, let $G_b(P)$  be $\sk(P)$ with all edges being bidirected. We shall need the following results from \cite{sad17}:

\begin{lemma}\label{prop:und}
Probability distribution $P$ and $G_u(P)$ are faithful if and only if $P$ satisfies intersection, singleton-transitivity, and upward-stability.
\end{lemma}
\begin{lemma}\label{prop:bid}
Probability distribution $P$ and $G_b(P)$ are faithful if and only if $P$ satisfies composition, singleton-transitivity, and downward-stability.
\end{lemma}

\subsubsection{Gaussian graphical models} Gaussian graphical models exemplify the properties describe above.
Let $P$ be a regular multivariate Gaussian distribution in $\mathbb{R}^V$ with covariance matrix $\Sigma$. It then holds that $\mathcal{J}(P)$ is a compositional graphoid. This follows from the fact that for such a distribution
\begin{equation}\label{eq:3}
A\ci B\cd C \iff (\Sigma_{ABC})^{-1}_{A,B}=0,
\end{equation}
where $(\Sigma_{ABC})^{-1}$ is the \emph{concentration matrix} of the distribution of $X_{A\cup B\cup C}$; for detailed proofs, see \citet{stu05}.

Define the undirected graph $G_u$ as follows: there is no edge between nodes $i$ and $j$ if $\Sigma^{-1}_{ij}=0$; and similarly define the bidirected graph $G_b$ as follows: no edge between $i$ and $j$ if $\Sigma_{ij}=0$. Hence, the pairwise Markov properties are satisfied with respect to  $G_u$ and $G_b$, and by the discussions above $P$ is Markov to these graphs.

\section{DPPs and conditional independence}
In this section, we first focus on the context-specific conditional independence models where the conditioning set $X_C$ is equal to the vector $1_C$, but eventually we provide results for general independence as well. First, we prove the following linear algebraic result:

\begin{lemma}\label{lem:4}
Let $K$ be a symmetric positive semidefinite matrix $K$ with row and column vector $V$. For disjoint, non-empty subsets $A$, $B$ and $C$ of $V$, the following statements are equivalent:
\begin{enumerate}
 \item $\det(K_{ABC})\det(K_C)=\det(K_{AC})\det(K_{BC})$;
 \item $(K_{ABC})^{-1}_{A,B}=0$;
 \item $K_{A,B} = K_{A,C}  K^{-1}_C K_{C,B} $.
\end{enumerate}
Simlarly, $\det(K_{AB})=\det(K_{A})\det(K_{B})$ if and only if $(K_{AB})_{A,B} = (K_{AB})^{-1}_{A,B} = 0$.
\end{lemma}
\begin{proof}
Let $K^C_{ABC}$ be the Schur complement of $K_{C}$ in $K_{ABC}$.
By the properties of the Schur complement, $(K_{ABC})^{-1}_{AB} = (K^C_{ABC})^{-1}$. Thus, $(K_{ABC})^{-1}_{A,B} = 0$ if and only if $(K^C_{ABC})^{-1}_{A,B} = 0$. Since $K^C_{ABC}$ is symmetric and its entries are indexed by $A \cup B$, we see that, if $(K^C_{ABC})^{-1}_{A,B} = 0$, then $(K^C_{ABC})^{-1}$ must be block diagonal. In turn, this is equivalent to  $(K^C_{ABC})_{A,B} = 0$.
By Theorem 15 in \cite{tad14}, this is then equivalent to the events $\{ X_A = 1_A\}$ and $\{ X_B = 1_B\}$ being conditionally independent given $\{ X_C = 1_C\}$, or, by definition, to the identity
\[
\frac{\det(K_{ABC})}{\det(K_C)} = \frac{ \det(K_{AC})}{\det(K_C)} \frac{\det(K_{BC})}{\det(K_C)},
\]
which simplifies to $\det(K_{ABC})\det(K_C)=\det(K_{AC})\det(K_{BC})$. Thus, {\it 1} and {\it 2} are equivalent.

To prove that  {\it 2} is equivalent to {\it 3} notice that,  since $K^C_{ABC} = K_{AB} - K_{AB,C} K_C^{-1}K_{C,AB}$, we have that $(K^C_{ABC})_{A,B} = K_{A,B} - K_{A,C}  K^{-1}_C K_{C,B}$. Thus $(K^C_{ABC})^{-1}_{A,B} = 0$ is equivalent to $K_{A,B} = K_{A,C}  K^{-1}_C K_{C,B} $.

The final claim follows from Corollary 7 in \cite{tad14} and the fact that $(K_{AB})_{A,B} =0$ if and only if $(K_{AB})^{-1}_{A,B} = 0$.
%
%
%
%
%
\end{proof}

Condition 2. in the previous result is analogous to the right clause in the equivalence \eqref{eq:3}, with the kernel $K$ playing the role of the covariance matrix $\Sigma$. This would suggest that, for a discrete DPP, condition 2. is equivalent to $X_A$ being conditionally independent to $X_B$ given $X_C$. This is in fact true provided that $X_C =1_C$, as we demonstrate next.

The proof of the following result uses the same idea as the proof of Theorem 1 in \citet{drt08}.
\begin{theorem}\label{thm:1}
For a DPP with kernel $K$ and any triplet $A$, $B$ and $C$ of non-empty disjoint subsets  of $V$, $(K_{ABC}^{-1})_{A,B}=0$ if and only if $X_A\ci X_B\cd X_{C}=1_C$.
\end{theorem}
\begin{proof}
By Lemma \ref{lem:4}, $(K_{ABC}^{-1})_{A,B}=0$ is equivalent to  $\det(K_{ABC})\det(K_C)=\det(K_{AC})\det(K_{BC})$.
The proof of Lemma \ref{lem:4}, in turn shows that this is equivalent to $P(X_A=1_A,X_B=1_B\cd X_C=1_C)=P(X_A=1_A\cd X_C=1_C)P(X_B=1_B\cd X_C=1_C)$.

The last identity will provide the basis for an inductive argument that will show that, if $(K_{ABC}^{-1})_{A,B}=0$, then $P(X_A=i_A,X_B=i_B\cd X_C=1_C)=P(X_A=i_A\cd X_C=1_C)P(X_B=i_B\cd X_C=1_C)$, for arbitrary $i_A \in \{0,1\}^A$ and $i_B \in \{0,1\}^B$. This is precisely the conditional probability clause $X_A\ci X_B\cd X_{C}=1_C$. The induction will be on the number of $0$'s in $(i_v)_{v\in A\cup B}$, with the base case -- where there are no $0$'s -- is proven above.



Thus, suppose that the result is true for less than $r$ $0$'s, and assume that there are $r$ $0$'s in $(i_v)_{v\in A\cup B}$. Suppose also that $i_v$=0. Without loss of generality assume that $v\in A$. We have that
\begin{equation*}
\begin{split}
P(X_A=i_A,X_B=i_B\cd X_C=1_C)\\=P(X_{A\setminus\{v\}}=i_{A\setminus\{v\}},X_B=i_B\cd X_C=1_C)\\-P(X_{A\setminus\{v\}}=i_{A\setminus\{v\}},X_B=i_B,X_v=1\cd X_C=1_C)\\=(P(X_{A\setminus\{v\}}=i_{A\setminus\{v\}}\cd X_C=1_C)\\-P(X_{A\setminus\{v\}}=i_{A\setminus\{v\}},X_v=1\cd X_C=1_C))\\P(X_B=i_B\cd X_C=1_C)\\=P(X_A=i_A\cd X_C=1_C)P(X_B=i_B\cd X_C=1_C),
\end{split}
\end{equation*}
where the second equality follows from the induction hypothesis. Thus, $(K_{ABC}^{-1})_{A,B}=0$  implies that $X_A\ci X_B\cd X_{C}=1_C$.

To show the reverse implication, assume that $X_A\ci X_B\cd X_{C}=1_C$. Then, the events $X_A = 1_A$ and $X_B = 1_B$ are conditionally independent given $X_C = 1_C$, which is equivalent to $\det(K_{ABC})\det(K_C)=\det(K_{AC})\det(K_{BC})$. The result now follows  from Lemma \ref{lem:4}.
\end{proof}


%

As an immediate corollary, we obtain the following context-specific pairwise Markov property.

\begin{coro}\label{coro:2}
For DPPs, $X_i\ci X_j\cd X_{R}=1_R$, where $R=V\setminus\{i,j\}$, if $(K^{-1})_{i,j}=0$.
\end{coro}

Using the fact that the complement of a DPP with kernel $K$ is a DPP with kernel $I-K$, we see that that the claims in both  Theorem \ref{thm:1} and Corollary \ref{coro:2} hold true if the conditioning events are $X_C = 0_C$ and $X_{R}=0_R$, respectively, provided $K$ is replaced by $I-K$.

Next, by using very similar and in fact simpler arguments as in the proof of Theorem \ref{thm:1} for the case of $C = \emptyset$, we arrive at important result characterizing marginal independence in discrete DPPs based on the zero entries in $K$. See also Theorem 8 in \cite{tad14}.

\begin{theorem}\label{coro:1}
For DPPs, $X_A\ci X_B$ if and only if $K_{A,B}=0$. In particular, $X_i\ci X_j$ if and only if $k_{ij}=0$.
\end{theorem}

Interestingly, since the correlation between any two variables $X_i$ and $X_j$ is
\[
 \frac{ - K^2_{ \{i,j\} } }{ \sqrt{ K_{ \{i\} } ( 1 - K_{\{i\}}) K_{ \{j\} }(1 - K_{\{j\}})}},
\]
Theorem \ref{coro:1} further implies that
\begin{coro}\label{cor:corr}
For a DPP, $X_i\ci X_j$ if and only if $X_i$ and $X_j$ are uncorrelated.
\end{coro}

Finally, we show that the context-specific conditional independence statements in Theorem \ref{thm:1} fulfill all the desiderata of a probabilistic model, much like the Gaussian graphical models.

\begin{prop}\label{prop:1}
For a DPP $P$, the context-specific independence model $\mathcal{J}_{\cdot\cd X_C=1_C}(P)$  satisfies intersection, composition, and singleton-transitivity.
\end{prop}
\begin{proof}
Because of Theorem \ref{thm:1} and Equation (\ref{eq:3}), the proof is the same as the proof for the independence model induced by Gaussian distribution satisfying these properties.
%
%
%
\end{proof}
Notice that because the complement of a DPP is a DPP, the above result also holds for the context-specific independence model $\mathcal{J}_{\cdot\cd X_C=0_C}(P)$.

In fact, except intersection, the other two properties are also satisfied for the general DPP induced independence models (In fact singleton-transitivity is always satisfied by binary distributions; see Example 3.1.7 of \cite{drts09}; but, here we provide a short proof for this special case):
\begin{prop}\label{prop:2}
For a DPP $P$, the independence model $\mathcal{J}(P)$  satisfies  composition and singleton-transitivity.
\end{prop}
\begin{proof}
In order to prove composition, first notice that since DPPs are closed under marginalization, we only need to focus on $X_{A\cup B\cup C\cup D}$. Also, since DPPs are closed under conditioning on any binary vector, we only need to check that $A\ci B$ and $A\ci D$ imply $A\ci B\cup D$. This follows from Theorem \ref{coro:1}.

To prove singleton-transitivity, with the same argument as for composition, we only need to prove that $i\ci j$ and $i\ci j\cd l$ imply $i\ci l$ or $j\ci l$.  We have, by Theorem \ref{coro:1}, that $k_{ij}=0$. If $l=1$ then, by Lemma \ref{lem:4} and the inverse of a $3\times 3$ matrix, $k_{ij}k_{ll}=k_{il}k_{jl}$. Therefore, either $k_{il}=0$ or $k_{jl}=0$, which implies the result. If $l=0$ then the result follows from the fact that the complement of a DPP is a DPP.
\end{proof}
However, notice that, for a DPP $P$, the independence model $\mathcal{J}(P)$ does not necessarily satisfy downward-stability. Examples for this are plentiful; for example, let
$$K=\left(
      \begin{array}{ccc}
         0.11& 0.04 & -0.10 \\
        0.04 & 0.29 & -0.22 \\
        -010 & -0.22 & 0.54 \\
      \end{array}
    \right)$$
This implies that $X_1\notci X_2$, but $K^{-1}_{1,2}\approx 0$, which implies that $X_1\ci X_2\cd X_3$.

%

\section{DPPs and Graphical Models}
We will now use the results from the previous section to show how discrete DPPs are Markov to bidirected graphs corresponding to the non-zero entries of $K$ and that context-specific discrete DPPs (conditioning on events of the form $X_C = 1_C$ or $ X_C = 0_C$) are Markov with respect to the undirected graph corresponding to the non-zero entries of $K^{-1}$.
\subsection{DPPs and Bidirected Graphical Models}
Given a DPP $P$ with kernel $K$, define $G_b(P)$ to be the bidirected graph with node set $V$ where there is no edge between nodes $i$ and $j$ if and only if $K_{ij}=0$ and vice versa.

By Theorem \ref{coro:1}, we immediately have that
\begin{prop}\label{prop:3}
A DPP $P$ is pairwise Markov to $G_b(P)$.
\end{prop}

In fact, we can obtain a stronger result.

\begin{prop}
A DPP $P$ is Markov to $G_b(P)$.
\end{prop}
\begin{proof}
The proof is a direct corollary of Lemma \ref{lem:3} below since (\ref{eq:2}) and consequently the connected set Markov property is satisfied.
\begin{lemma}\label{lem:3}
Let $P$ be a DPP with kernel $K$. For every subgraph $D$ of $G_b(P)$, it holds that $q_D=q_{C_1}q_{C_2}\cdots q_{C_r}$, where $C_1,\dots,C_r$ are the maximal connected components of $D$.
\end{lemma}

The proof of the previous claim follows from the fact that $K_D$ is a block diagonal matrix with blocks $C_1,\dots,C_r$, and hence, $\det(K_D)=\det(C_1)\dots\det(C_r)$.

\end{proof}

This global Markov property can also be proven  using the pairwise Markov property (Proposition \ref{prop:3}) and the composition property (Proposition \ref{prop:2}).

Finally, making use of Lemma \ref{prop:bid}, we have the following result regarding faithfulness:
\begin{prop}
A DPP $P$ is faithful to $G_b(P)$ if and only if $\mathcal{J}(P)$ satisfies downward-stability.
\end{prop}
As the example at the end of the previous section shows, downward-stability does not always hold, and, although DPPs are negatively associated,  faithfulness does not immediately follow from Markovness. This result should be contrasted with the conditions for faithfulness in $\mathrm{MTP}_2$ distributions, which are instead positively associated; see Section 6 in \cite{fal15}.

\subsection{DPPs and Undirected Graphical Models}
Given a DPP $P$ with kernel $K$, define $G_u(P)$ to be the undirected graph on $V$ where there is no edge between nodes $i$ and $j$ if and only if $K^{-1}_{ij}=0$.{}
By Corollary \ref{coro:2}, we have the following:
\begin{prop}
For a DPP $P$, the context-specific independence model $\mathcal{J}_{\cdot\cd X_C=1_C}(P)$ is pairwise Markov to $G_u(P)$.
\end{prop}

Like in the unconditional case, we caen easily strangthen the previous resul to obtain global Markov properties.

\begin{prop}
For a DPP $P$, the context-specific independence model $\mathcal{J}_{\cdot\cd X_C=1_C}(P)$ is Markov to $G_u(P)$.
\end{prop}
\begin{proof}
The proof follows from the fact that $\mathcal{J}_{\cdot\cd X_C=1_C}(P)$ satisfies symmetry and decomposition (Lemma \ref{lem:2}) as well as the intersection property (Proposition \ref{prop:1}).
\end{proof}
By the use of Lemma \ref{prop:und}, we also have the following:
\begin{prop}
For a DPP $P$, the context-specific independence model $\mathcal{J}_{\cdot\cd X_C=1_C}(P)$ is faithful to $G_u(P)$ if and only if $\mathcal{J}(P)$ satisfies upward-stability.
\end{prop}
In a similar fashion to downward-stability, it is easy to find examples of DPPs whose induced context-specific independence model does not satisfy upward-stability.

\section{Summary}
We have shown that the bidirected DPP graphical model behaves similarly to the bidirected Gaussian graphical model. To us, this is  noteworthy, as the only other known binary graphical model with the desirable properties of the Gaussian graphical model is the very restrictive class of symmetric binary models \citep{wer09}, where the marginal probabilities of $0$ and $1$ are equal.
We have also illustrated that, for undirected graphs, it is the context-specific induced independence model by DPP given the vector of $1$'s (or $0$'s), which acts similarly to the independence models induced by the Gaussian distribution. This provides a nice example for the use of context-specific independence models as opposed to the general ones.

Although DPPs are negatively associated \citep{lyo03}, the faithfulness property in the bidirected graph case is not satisfied. Hence, DPPs are not (in some sense) dual to positively associated $\mtp$ distributions \citep{Kar80} (which are always faithful to their corresponding undirected graphs \citep{fal15}). It would be interesting to find sufficient (and even necessary) conditions on the kernel $K$ so that downward-stability is satisfied. Such a result would imply  that the restricted set of kernels based on these conditions generates DPPs that are faithful to their corresponding bidirected graph. Similarly, it would be interesting to find conditions on $K$ so that the upward-stability is satisfied. This implies faithfulness of the context-specific independence model and the corresponding undirected graph.

In this note we only discuss undirected and bidirected graphs associated to DPPs, but the fact that context-specific  independence models induced by discrete DPPs have essentially all the independence properties of the Gaussian distribution implies that graphical models arising from other types of graphs -- such as directed acyclic graph or mixed graph \citep[see][]{sadl16} -- can be proposed for DPPs. The latter could cover the case where there are latent variables in a DPP graphical model.
\subsubsection*{Acknowledgements}
The authors are grateful to Shaowei Lin for providing the idea for proving Lemma \ref{lem:4}.

\bibliography{bib}
\end{document}